\documentclass [10pt,twoside,b5paper]{article}
\usepackage{amsfonts}
\usepackage{amsthm}
\usepackage{amsmath}
\usepackage{amstext}
\usepackage{amssymb}
\usepackage{mathrsfs}
\usepackage{epsf}              
\usepackage{graphicx}          
\usepackage{fancybox}          
\usepackage{color}             
\usepackage{fancyhdr}
\usepackage[hang,footnotesize]{caption2}  
\def\mathcal{\mathscr}
\newfont{\aaa}{cmb10 at 18pt}
\newfont{\bbb}{cmb10 at 10pt}

\newtheorem{theo}{Theorem}
\newtheorem{defi}[theo]{Definition}

\def\n{\nabla}

\def\pt{\partial}
\def\n{\mathbb{N}}
\def\z{\mathbb{Z}}
\def\c{\mathbb{C}}
\def\dim{\hbox{dim}}
\def\ad{\hbox{ad}}

\def\F{\mathbb F}

\def\a{\alpha}
\def\b{\beta}
\def\sg{\sigma}
\def\ll{\lambda}

\newfont{\df}{eufm10}
\def\vep{\varepsilon}

\def\sg{\sigma}

\def\ll{\lambda}
\def\der{\hbox{Der}\,}
\def\aut{\hbox{Aut}\,}

\def\vp{\varphi}

\def\r{\gamma}

\def\de{\delta}
\def\dim{\hbox{\rm dim}\,}

\def\hom{\hbox{\rm Hom}\,}

\def\span{\hbox{\rm Span}\,}
\def\ad{\hbox{\rm ad}\,}

\def\mg{{\bf \frak g}}

\newcommand{\beq}{\begin{equation}}
\newcommand{\eeq}{\end{equation}}
\newcommand{\bey}{\begin{eqnarray}}
\newcommand{\eey}{\end{eqnarray}}
\newcommand{\beyy}{\begin{eqnarray*}}
\newcommand{\eeyy}{\end{eqnarray*}}

\renewcommand{\headsep}{0.7cm}
 \makeatletter
\def\@evenhead{
   \vbox{\hbox to \textwidth
{}{\hspace{0mm}{\footnotesize \thepage}}{\hspace{8cm}
   {\footnotesize {Dong Liu, Linsheng Zhu}}}
  \protect\vspace{1truemm}\relax
   \hrule depth0pt height0.15truemm width\textwidth
   }}
   \def\@evenfoot{}
\def\@oddhead{
    \vbox{\hbox to \textwidth
  {{\hspace{0cm}{\footnotesize  The generalized Heisenberg-Virasoro algebra} \hfill{\footnotesize \thepage}}\hspace{0mm}}{}
   \protect\vspace{1truemm}\relax
   \hrule depth0pt height0.15truemm width\textwidth
  }}
  \def\@oddfoot{}
\makeatother

\begin{document}

\thispagestyle{empty} \thispagestyle{fancy} {
  \fancyhead[lO,RE]{\footnotesize  Front.  Math. China \\
  DOI 10.1007/s11464-009-00-\\[3mm]
  \includegraphics[0,-50][0,0]{11.bmp}}
  \fancyhead[RO,LE]{\scriptsize \bf 
} \fancyfoot[CE,CO]{}}
\renewcommand{\headrulewidth}{0pt}
\renewcommand{\headsep}{0.7cm}


\setcounter{page}{1}
\qquad\\[8mm]

\noindent{\aaa{The Generalized Heisenberg-Virasoro algebra$^{^{^{\displaystyle*}}}$}}\\[1mm]

\noindent{\bbb Dong Liu$^{1},$\quad Linsheng Zhu$^2$}\\[-1mm]

\noindent\footnotesize{1\ \ Department of Mathematics, Huzhou Teachers
College,  Zhejiang Huzhou, 313000, China\\
2\ \ Department of Mathematics, Changshu Institute of
Technology, Jiangsu Changshu, 215500, China}\\[6mm]

 \vskip-2mm
 \noindent{\footnotesize$\copyright$ Higher Education Press and
Springer-Verlag 2009}
 \vskip 4mm

\normalsize\noindent{\bbb Abstract}\quad  In this paper, we mainly study the generalized Heisenberg-Virasoro
algebra. Some structural properties of
the Lie algebra are obtained.\vspace{0.3cm}

\noindent{\bbb Keywords}\quad The generalized Heisenberg-Virasoro
algebra,  central extension, automorphisms\\
{\bbb MSC}\quad 17B56; 17B68 \\[0.4cm]

\noindent{\bbb{1\quad Introduction}}\\[0.1cm]

\noindent  The twisted Heisenberg-Virasoro algebra $H_{Vir}$ has been first
studied by Arbarello et al. in Ref. [1], where a connection is
established between the second cohomology of certain moduli spaces
of curves and the second cohomology of the Lie algebra of
differential operators of order at most one:
$$L_{HV}=\{f(t)\frac{d}{dt}+g(t)|f,g\in{\mathbb
C}[t,t^{-1}]\}.$$

As a vector space over $\c$, $H_{Vir}$ has a basis $\{L(m), I(m),
C_L, C_I, C_{LI}, m\in\z\}$, subject to the following relations:
\begin{eqnarray*}
&&[L(m), L(n)]=(n-m)L(m+n)+\de_{m+n, 0}{1\over 12}(m^3-m)C_L;\\
&&[I(m), I(n)]=n\de_{m+n, 0}C_{I};\\
&&[L(m), I(n)]=nI(m+n)+\de_{m+n, 0}(m^2-m)C_{LI};\\
&&[H_{Vir}, C_L]=[H_{Vir}, C_I]=[H_{Vir}, C_{LI}]=0.
\end{eqnarray*}

Clearly the Heisenberg algebra $H=\c\{I(m), C_I\mid m\in\z\}$ and
the Virasoro algebra $Vir=\c\{L(m), C_L\mid m\in\z\}$ are
subalgebras of $H_{Vir}$.

Arbarello et al.
(in Ref. [1]) also proved that any irreducible highest weight module
for $H_{Vir}$ is isomorphic to the tensor product of an irreducible

\vspace{-2mm}
\noindent \hrulefill\hspace{117mm}\\
{\footnotesize $^*$  Received Oct. 10, 2008; accepted Feb. 25,
2009}

\vspace{-1mm}
\noindent \hrulefill\hspace{117mm}\\
{\footnotesize Corresponding author: Dong Liu, E-mail:
liudong@hutc.zj.cn}

\newpage

\ noindent  module for the Virasoro algebra and an irreducible module for the
infinite-dimensional Heisenberg algebra when the central element of
the Heisenberg subalgebra acts in a non-zero way. The structure of
the irreducible representations for $H_{Vir}$ at level zero was
studied in Ref. [3]. The Harish-Chandra modules over $H_{Vir}$ were
classified in Refs. [11, 12]. Some structure properties and
representations of the twisted Heisenberg-Virasoro Lie algebra were
obtained in Refs. [7, 8, 9, etc.].

Recently, a number of new classes of infinite-dimensional simple
Lie algebras over a field of characteristic 0 were discovered by
several authors. Among those algebras, are the generalized Witt
algebras, the generalized Virasoro algebras introduced in Ref. [16],
which are one-dimensional universal central extensions of some
generalized Witt algebras, and the Lie algebras of generalized
Weyl type introduced and studied in Refs. [15,17, etc.].

Motivated by the above algebras, we introduce a new Lie algebra, the
generalized Heisenberg-Virasoro algebra (see Definition 1 in Section 2), which
is a generalization of the twisted Heisenberg-Virasoro Lie algebra
$H_{Vir}$ from the integer ring $\z$ to an additive subgroup $M$ of
a field $\F$ (The special case $M=\z^2$ was studied in Ref. [18]). Its Verma modules were studied in  Ref. [14].
In this paper, we mainly study some properties of this Lie algebra.

The paper is organized as follows. In Section 2, we recall some
notions of generalized Witt algebras and the twisted
Heisenberg-Virasoro algebra, and then introduce the definition of
the generalized Heisenberg-Virasoro algebra. In Section 3, we prove
that the generalized Heisenberg-Virasoro algebra is the universal
extension of the Lie algebra of generalized differential operators
on a circle of order at least one: $L={\mathbb F}\{t^x\pt,\, t^y\mid
x, y\in A, \pt\in T\}$. In Section 4 and Section 5, we determine
derivations and automorphisms of the generalized Heisenberg-Virasoro
algebra.

Throughout this paper, ${\mathbb F}$ denotes a field of
characteristic zero, and $\F^*=\F\backslash\{0\}$. $\c$, $\z$,
$\n$, $\z^+$ denote the filed of complex numbers, the set of all
integers, the set of all nonnegative integers, the set of all
positive integers, respectively.
\\[4mm]

\noindent{\bbb 2\quad The generalized Heisenberg-Virasoro algebra}\\[1mm]

\noindent First we recall some notions of the generalized Witt algebras as
defined in Refs. [4, 10].

Let $A \,(\ne0)$ be an abelian group and $T$ a vector space over
${\mathbb F}$. In this paper, we are only interested in the case of
$\dim T=1$, so we always assume that $T={\mathbb F}\pt$ throughout
this paper. We denote by $\mathbb FA$ the group algebra of $A$ over
${\mathbb F}$. The elements $t^x, x\in A$, form a basis of the
algebra, and the multiplication is defined by $t^xt^y=t^{x+y}$. We
shall write 1 instead of $t^0$. The tensor product $W=\F A\otimes_\F
T$ is free left  $\F A$-module. We shall usually write $t^x\pt$
instead of $t^x\otimes\pt$. We now fix a pairing $\vp: T\times A\to
\F$, which is ${\mathbb F}$-linear in the first variable and
additive in the second one. For convenience we shall use the
following notations:
$$\vp(\pt, x)=\langle\pt, x\rangle=\pt(x)$$ for arbitrary $x\in A$.

There is a unique ${\mathbb F}$-linear map $W\times W\to W$
sending $(t^x\pt, t^y\pt)$ to
$$[t^x\pt, t^y\pt]:=t^{x+y}\pt(y-x)\pt \eqno(2.1)$$  for arbitrary $x, y\in A$.
This map makes  $W=W(A, T, \vp)$ into a Lie algebra which is called {\sl generalized Witt algebra}.
Setting $W_x=t^xT, x\in A$, then
$$W=\oplus_{x\in A}W_x$$ is a graded Lie algebra.

We say that $\vp$ is nondegenerate if $A_0:=\{x\in A: \pt(x)=0\}=0$.

The following theorem is due to N. Kawamoto (Ref. [10]).

\vspace{1mm}
\noindent{\bf Theorem 2.1} (Ref. [10])  {\sl  Suppose that characteristic of ${\mathbb F}$ is 0. Then
$W=W(A, T, \vp)$ is a simple Lie algebra if and only if $\vp$ is
nondegenerate}.

In this paper, we always suppose that $\vp$ is nondegenerate.

In Ref. [5], the second cohomology group of the generalized Witt
algebra was constructed.

\vspace{1mm}
\noindent{\bf Theorem 2.2} (Ref. [5]) {\sl Let $W=W(A, T, \vp)$ be a simple
generalized Witt algebra and $T=\F\pt$ one-dimensional, then
$H^2(W, \F)$ is $1$-dimensional and is spanned by the cohomology
class $[\psi]$, where $\psi: W\times W\to \F$ is the $2$-cocycle
defined by}
$$\psi(t^x\pt, t^y\pt)=\delta_{x+y, 0}\pt(x)^3, \ x, y\in A. \eqno(2.2)$$

The generalized Virasoro algebra is the universal central
extension of the generalized Witt algebra when $\vp$ is
nondegenerate (Ref. [16]).

By definition, the {\it generalized Virasoro algebra} $Vir$ is a
Lie algebra generated by $\{L(x)=t^x\pt, \, x\in A$\} and $C_{L}$
subject to the following relations:
\begin{eqnarray*}
&&[L(x), L(y)]=\pt(y-x) L(x+y)+\de_{x+y, 0}{1\over12}(\pt(x)^3-\pt(x))C_L\\
&&[L(x), C_L]  = 0.
\end{eqnarray*}

Now we introduce the algebra of generalized differential operators
in Refs. [15, 16].

Let ${\mathcal D}$ be the associative algebra generated by the
elements $\{t^x, x\in A\}$ and the element $\pt$, subject to the
following relation:
$$\pt t^x-t^x\pt=\pt(x)t^x,  \quad \forall x\in A. \eqno(2.3)$$
The algebra ${\mathcal D}$ is called the {\it algebra of
generalized differential operators}. Denoted by ${\mathcal D}^-$
the Lie algebra associated with ${\mathcal D}$.

\vspace{1mm}
\noindent{\bf Example.} Let $A=\z$, then the group algebra $\F A$
becomes identified with the algebra of Laurent polynomials $\F[t,
t^{-1}]$. Let $T$ be a vector space with basis $\{\pt\}$. Define
the pairing $\vp: T\times A\to \F$ by setting $\langle\pt,
m\rangle=m$. Then $\vp$ is nondegenerate and the algebra of
generalized differential operators is just the algebra of
differential operators with coefficients in the Laurent polynomial
ring.

As a vector space over ${\mathbb F}$, ${\mathcal D}$ has a basis
$\{t^x\pt^m, x\in A,\, \pt\in T, m\in\n\}$. It is easy to see that
the following relation holds in the associative algebra ${\mathcal
D}$.
$$(t^x\pt^m)(t^y\pt^n)=t^{x+y}\sum_{i=0}^m {m\choose i}\pt(y)^i\pt^{m+n-i},  \quad \forall x, y\in A,  m, n\in\n. \eqno(2.5)$$

By (2.5) we can see that the generalized Witt algebra $W=W(A, T,
\vp)$ is a Lie subalgebra of ${\mathcal D}^-$, the Lie algebra of
generalized differential operators.

The above results lead us to the following definition.
\begin{defi}
Let $T={\mathbb F}\partial$ and $W=W(A,T,\varphi)$ be a simple
generalized Witt algebra. The generalized Heisenberg-Virasoro
algebra ${\mathcal L}$ is a Lie algebra generated by
$\{L(x)=t^x\pt, I(x)=t^x, C_L, C_I, C_{LI}\}$, subject to the
following relations:
\begin{eqnarray*}
&&[L(x), L(y)]=\pt(y-x)L(x+y)+\de_{x+y, 0}{1\over 12}(\pt(x)^3-\pt(x))C_L;\\
&&[I(x), I(y)]=\pt(y)\de_{x+y, 0}C_{I};\\
&&[L(x), I(y)]=\pt(y)I(x+y)+\de_{x+y, 0}(\pt(x)^2-\pt(x))C_{L};\\
&&[{\cal L}, C_L]=[{\cal L}, C_I]=[{\cal L}, C_{LI}]=0.
\end{eqnarray*}
\end{defi}

The Lie algebra ${\mathcal L}$ has a generalized Heisenberg
subalgebra and a generalized Virasoro subalgebra interwined with a
2-cocycle.

Moreover, we shall prove that the generalized Heisenberg-Virasoro
algebra ${\mathcal L}$ is the universal central extension of the Lie
algebra of generalized differential operators of order at least one:
${\cal D}_1={\mathbb F}\{t^x\pt,\,  t^y\mid x, y\in A,\, \pt\in
T\}$.
\\[4mm]
\noindent {\bbb 3\quad The universal central extension of the generalized Heisenberg-Virasoro
algebra}\\[1mm]

\noindent Now we consider the Lie subalgebra ${\cal D}_1=\mathbb F\{t^x\pt,\,
t^y\mid x, y\in A\}$ of ${\mathcal D}^-$, the Lie algebra of
generalized differential operators.

\vspace{1mm}
\noindent{\bf Theorem 3.1} $\dim H^2({\cal D}_1, \F)=3.$

\noindent{\bf Proof.} Let $\a$ be any 2-cocycle of ${\cal D}_1$ and
we set $\a_{x, y}=\a(t^x\pt, t^y\pt)$, $\a_x^y=\a(t^x\pt, t^y)$ and
$\a^{x, y}=\a(t^x, t^y)$. Considering $\a([t^z\pt, t^x], t^y)$,
$\a([t^z\pt, t^x\pt], t^y)$, $\a([t^z\pt, t^x\pt], t^y\pt)$ and
using the relations:
$$[t^x\pt, t^y\pt]=\pt(y-x)t^{x+y}\pt, \quad [t^x\pt, t^y]=\pt(y)t^{x+y}, \quad [t^x, t^y]=0,$$
and the cocycle rule, we get the following
$$\pt(x)\a^{z+x, y}=\pt(y)\a^{z+y, x},\eqno(3.1)$$
$$\pt(x-z)\a_{z+x}^y=\pt(y)(\a_z^{x+y}-\a_x^{y+z}),\eqno(3.2)$$
$$\pt(x-z)\a_{z+x, y}+\pt(y-x)\a_{x+y, z}+\pt(z-y)\a_{y+z, x}=0.\eqno(3.3)$$
The relation (3.1) immediately implies
$$\a^{x, y}={\pt(y)\over\pt(x_0)}\de_{x+y, 0}\a^{-x_0, x_0}, \forall x, y\in
A,\eqno(3.4)$$ for some $x_0\in A$ such that $\pt(x_0)\ne 0$.

By (3.4), $\a^{x, y}$ is determined by nontrivial 2-cocycle $\psi_1:
{\cal D}_1\times {\cal D}_1\to \F$ defined by $\psi_1(t^x,
t^y)=\de_{x+y, 0}\pt(y)$, for all $x, y\in A$.

From the proof of Theorem 6.1 in Ref. [5], we see that, up to a
co-boundary, it follows that $\a_{x, y}$ is determined by nontrivial
2-cocycle $\psi_2$, where $\psi_2(t^x\pt, t^y\pt)=\de_{x+y,
0}(\pt(x)^3-\pt(x)),$ for all $x, y\in A$.

Setting $y=0$ in (3.2) we deduce that
$$\a_x^0=0,\quad \forall x\in A. \eqno(3.5)$$

Let $g: {\cal D}_1\to \F$ be the linear function $g(t^x\pt)=0$ and
$g(t^x)={1\over \pt(x)}\a(\pt, t^x)$ if $x\ne 0$ and $g(t^0)=0$.
Then $g$ induces a 2-coboundary $\psi_g(u, v)=g([u, v])$, $u, v\in
{\cal D}_1$. Hence by replacing $\a$ with the 2-coboundary
$\a-\psi_g$, we may assume that $\a(\pt, t^y)=0$ for all $y\ne
0$. Hence by (3.5), we obtain
$$\a_0^x=0,\quad \forall x\in A.$$
Setting $z=0$ in (3.2), we have
$$\a_x^y=0, \quad \forall x,y\in A\ \hbox{and}\ x+y\ne 0.\eqno(3.6)$$

On the other hand, setting $y=-z-x$ in (3.2), we have
$$\pt(x-z)\a_{x+z}^{-x-z}=\pt(x+z)(\a_x^{-x}-\a_z^{-z}), \quad \forall x, z\in A. \eqno(3.7)$$

Let $x\in A$ be such that $\pt(x)\ne0$. By substituting $kx$ for $x$
and $lx$ for $z$, we obtain
$$(k-l)f(k+l)=(k+l)(f(k)-f(l)), \eqno(3.8)$$
where $f(k)=\a_{kx}^{-kx}$, $k\in\z$.

It follows from (3.8) that all the values $f(k)$ can be computed if
$f(1)$ and $f(2)$ are known. Hence the general solution (3.8) is
given by $f(k)=ak+bk^2$ for some constants $a, b\in \F$.

Setting $k=1$ and $k=2$, we have
$$\a_{kx}^{-kx}={k\over 2}(4\a_x^{-x}-\a_{2x}^{-2x})+{k^2\over 2}(\a_{2x}^{-2x}-2\a_x^{-x}) \eqno(3.9)$$
for all $k\in\z$, provided that $\pt(x)\ne0$.

Let $x_0$ be a nonzero element of $A$. Using similar considerations
in Page 661 of Ref. [5], we obtain
$$\a_{x}^{-x}={\pt(x)\over 2\pt(x_0)}(4\a_{x_0}^{-x_0}-\a_{2x_0}^{-2x_0})+{\pt^2(x)\over 2\pt^2(x_0)}(\a_{2x_0}^{-2x_0}-2\a_{x_0}^{-x_0}), \quad \forall x\in A. \eqno(3.10)$$

Then $\a_x^y$ is determined by a nontrivial 2-cocycle $\psi_3: {\cal
D}_1\times {\cal D}_1\to \F$ defined by $\psi_3(t^x\pt,
t^y)=\de_{x+y, 0}\pt(x)^2$ and the 2-coboundary $\psi_3'(t^x\pt,
t^y)=\de_{x+y, 0}\pt(x)$ for all $x, y\in A$. \hfill
$\rule[-.23ex]{1.0ex}{2.0ex}$

\vspace{1mm}
\noindent{\bf Corollary 3.2}
{\sl  The generalized Heisenberg-Virasoro algebra ${\mathcal L}$ is the
universal central extension of the Lie algebra ${\cal D}_1$ of
generalized differential operators of order at least one.}  \hfill
$\rule[-.23ex]{1.0ex}{2.0ex}$
\\[0.4cm]
\noindent{\bbb{4\quad Derivations of the generalized Heisenberg-Virasoro algebra}}\\[0.1cm]
\vspace{1mm}
\noindent   In this section we shall determine all derivations of the
generalized Heisenberg-Virasoro algebra ${\mathcal L}$.

First, we recall a result about derivations in Ref. [2].

\noindent{\bf Proposition 4.1} (Ref. [2]) {\sl
If $L$ is a perfect Lie algebra and $\hat L$ is a universal central extension of $L$,
then every derivation of $L$
lifts to a derivation of $\hat L$. If $L$ is centerless, the lift
is unique and} $\der(\hat L)= \der(L)$.

Due to Corollary 3.2 in Section 3 and Proposition 4.1, we
shall just determine all derivations of ${\cal D}_1={\mathbb
F}\{t^x\pt,\, t^y\mid x, y\in A, \pt\in T\}$, the Lie subalgebra
of ${\mathcal D}^-$.

We now describe two kinds of derivations of degree 0, which are
outer derivations of ${\cal D}_1$.

The linear maps $\sg_1, \sg_2,\sg_3: {\cal D}_1\to {\cal D}_1$
defined by
$$\sg_1(L(x))=\pt(x)I(x),\quad  \sg_1(I(x))=0, \quad \forall x\in
A,$$ and
$$\sg_2(L(x))=I(x),\quad  \sg_2(I(x))=0, \quad \forall x\in
A,$$ and
$$\sg_3(L(x))=0,\quad  \sg_3(I(x))=I(x), \quad \forall x\in
A,$$ are outer derivations of ${\cal D}_1$.

Let $\mu: A\to \F$ be an additive map, then the linear map
$\xi_{\mu}: {\cal D}_1\to {\cal D}_1$ defined by
$$\xi_{\mu}(L(x))=\mu(x)L(x),\quad \xi_{\mu}(I(x))=\mu(x)I(x), \quad \forall x\in A,$$
is also a derivation of degree 0.  It is clear that it is an outer
derivation if $\mu\ne k\pt$ for any $k\in\F$. In fact, if $\mu=k\pt$
for some $k\in\F$, then $\xi_\mu=\ad(k\pt)$.

Let $A$ be an abelian group, $\mg$ an $A$-graded Lie algebra, and
$V$ an $A$-graded left $\mg$-module. We also denote by $L_x$ resp.
$V_x$ the homogeneous components of ${\cal D}_1$ resp. $V$ for any
$x\in A$.

\vspace{1mm}
\noindent{\bf Proposition 4.2} (Refs. [5, 6])  {\sl  Every derivation $D\in \der(\mg, V)$ can be
written as $$D=\sum_{x\in A}D_x, \quad D_x\in \der(\mg, V)_x,$$ in
the sense that for every $v\in \mg$ only finitely many $D_xv\ne 0$
and $$Dv=\sum_{x\in A}D_xv,$$ where $\der(\mg, V)_x=\{\sg\in
\der(\mg, V)\mid \sg(L_y)\subset V_{x+y}, \forall y\in A\}$}.

A derivation $D: \mg\to V$ is called {\it locally inner}
(Ref. [5]) if it is a sum (may be infinite sum) of inner
derivations.

\vspace{1mm}
\noindent{\bf Proposition 4.3}  (Refs. [5, 6]) {\sl Suppose that the following conditions hold.

(1) $H^1(\mg_0, V_x)=0$ for $x\ne 0$;

(2) $\hom_{\mg_0}(L_x, V_y)=0$ for $x\ne y$.

Then $\der(\mg, V)_x$, $x\ne0 $, consists of inner derivations and
consequently
$$\der(\mg, V)=\der(\mg, V)_0+\der'(\mg, V),$$
where $\der'(\mg, V)$ is the space of locally inner derivations
$\mg\to V$.}

We  now use Proposition 4.3 to prove the following result.

\vspace{1mm}
\noindent{\bf Proposition 4.4}
$$\der({\cal D}_1)=\ad({\cal D}_1)+\der({\cal D}_1)_0.$$

\noindent{\bf Proof.}  By direct calculation, see Ref. [5].
\hfill  $\rule[-.23ex]{1.0ex}{2.0ex}$

\vspace{1mm}
\noindent{\bf Theorem 4.5} {\sl
Assume that $A\ne 0$ and $\vp$ is nondegenerate.  Then}
$$\der({\cal D}_1)_0=\ad({\cal D}_1)_0+\sum_{\mu\in\hom(A, \F)}\F\xi_\mu+\sum_{i=1}^3\F\sg_i,$$

\noindent{\bf Proof.}  Since $\vp$ is nondegenerate, we can
suppose that $\pt(x)\ne 0$ if $x\ne0$.

For any $D\in \der({\cal D}_1)_0$ and $x\in A $, we suppose
$$D(t^x)=t^x(\a_x\pt+\b_x), \a_x, \b_x\in{\mathbb F}. \eqno(4.1)$$
For any $y\in A$, by the fact that $[t^x, t^y]=0$, we obtain that
$$\a_x\pt(y)-\pt(x)\a_y=0.\eqno(4.2)$$
Choose $x_0\in A$ such that $\pt(x_0)=1$ and set $a_0=\a_{x_0}$. So
by (4.2) we have
$$\a_y=\pt(y)a_0,\eqno(4.3)$$
and (4.1) becomes
$$D(t^x)=t^x(a_0\pt(x)\pt+\b_x),\  a_0,\, \b_x\in{\mathbb F}.\eqno(4.4)$$
Now we suppose that
$$D(t^x\pt)=t^x(\r_x\pt+\ll_x), \r_x, \ll_x\in{\mathbb F}. \eqno(4.5)$$
Applying $D$ to $[t^y\pt, t^x]=\pt(x)t^{x+y}$, we obtain
$$[D(t^y\pt), t^x]+[t^y\pt, D(t^x)]=\pt(x)D(t^{x+y}),$$  i.e.
$$\r_y\pt(x)t^{x+y}+\b_x\pt(x)t^{x+y}+\pt(x)\pt(x-y)a_0t^{x+y}\pt=\pt(x)\pt(x+y)a_0t^{x+y}\pt+\pt(x)\b_{x+y}t^{x+y}.$$
So $\pt(x)\pt(x-y)a_0=\pt(x)\pt(x+y)a_0$ and
$\r_y\pt(x)+\b_x\pt(x)=\pt(x)\b_{x+y}$ for all $x, y\in A$. Hence
$$a_0=0\eqno(4.6)$$ and $$\b_{x+y}=\b_x+\r_y, \  \forall x\ne0, y\in A. \eqno(4.7)$$

Substitute $D$ by $D-\b_0\sg_3$, we can suppose $\b_0=0$. Setting
$y=-x$ in (4.7), we have $$\b_x=-\r_{-x}, \ \forall  x\in A.
\eqno(4.8)$$

So (4.4) and (4.5) become
$$D(t^x)=\b_xt^x.\eqno(4.9)$$
$$D(t^x\pt)=t^x(\r_x\pt+\ll_x),\ \b_x=-\r_{-x}, \ll_x\in{\mathbb F}. \eqno(4.10)$$

Applying $D$ to $[t^x\pt, t^y\pt]=\pt(y-x)t^{x+y}\pt$, we obtain
$$\r_{x+y}=\r_x+\r_y, \quad \hbox{and}\quad  \pt(y-x)\ll_{x+y}=\pt(y)\ll_y-\pt(x)\ll_x.\eqno(4.11)$$
By (4.11)and (4.8) we have $\r_x=\b_x$ since $\r_0=0$. Then $\b: A\to \F$ defined by $ \b(x)=\b_x$ is an
additive map. Moreover by using the method as in Section 3 (see (3.7)-(3.9))
 we can deduce that $\ll_x=a\pt(x)+b$ for some $a, b\in\F$ from (4.11). Replacing
$D$ by $D-\xi_{\b}-a\sg_1-b\sg_2$, we infer that $D=0$. \hfill
$\rule[-.23ex]{1.0ex}{2.0ex}$
\\[4mm]

\noindent{\bbb{5\quad The automorphism group of ${\mathcal L}$}}\\[0.1cm]

\noindent   Now, we recall a result about automorphisms in Ref. [13].

\vskip1mm
\noindent{\bf Proposition 5.1}  (Ref. [13]) {\sl
Let $L$ be a perfect Lie algebra and $\hat L$ be its universal
central extension. Every automorphism $\theta$ of $L$ admits a
unique extension $\hat\theta$ to an automorphism $\hat L$.
Furthermore, the map $\theta\to\hat\theta$ is a group
monomorphism}.

Due to Corollary 3.2 in Section 3 and Proposition 5.1,
we shall just determine the automorphism group of ${\cal
D}_1={\mathbb F}\{t^x\pt,\ t^y\mid x, y\in A,\, \pt\in T\}$. Set
$I={\mathbb F}\{t^x\mid x\in A\}$, then $I$ is the
unique maximal solvable ideal of ${\cal D}_1$. Clearly, if $\pi$
is an automorphism of ${{\cal D}_1}$, then
$$\pi(I)=I. \eqno(5.1)$$ Denote by ${\frak r}$ the set of all
inner automorphisms of ${\cal D}_1$, then ${\frak r}$ is a normal
subgroup of $\aut({\cal D}_1)$ and ${\frak r}$ is generated by
exp$(k\,\ad t^x),\ x\in A$, $k\in{\mathbb F}$.

For convenience, denote by $X(A)$ the group of characters of $A$,
i.e., the group homomorphisms $A\to {\mathbb F}^*$. Set ${\mathcal
E}=\{\vep\in\mathbb F^*\mid \span_{\mathbb F}\{t^{\vep x}, x\in
A\}=\span_{\mathbb F}\{t^{x}, x\in A\}\}$.

\vskip1mm
\noindent{\bf Remark.}  For example, if $A=\z$, then ${\cal
E}=\{1, -1\}$ (see Ref. [14]); if $A=\mathbb Q$, then ${\cal
E}=\mathbb Q^*$.

\vskip1mm
\noindent{\bf Lemma 5.2} {\sl
 Let $\pi\in
\aut({{\cal D}_1})$, then there exists $\eta\in {\frak r}$, such
that
$$\eta^{-1}\pi(\pt)=\vep^{-1}\pt+a,$$ for some $\vep\in {\cal E}$ and $a\in \mathbb F$}.

\noindent{\bf Proof.} Assume that
$$\pi(\pt)=\sum\limits_{z\in A_1}\lambda_zt^z\pt+
\sum\limits_{w\in A_2}\gamma_wt^w,$$ where $\lambda_z, \gamma_w\in
\mathbb{F}$ and $A_1, A_2$ are finite subsets of $A$. By (5.1), we
can assume
$$\pi(t^x)=\sum\limits_{y\in A_3}\nu(y)t^y,$$ where $A_3$ is a finite subset of
$A$.
Since $[\pt,t^x]=\pt(x)t^x,$ we have
$$[\sum\limits_{z\in A_1}\lambda_zt^z\pt+
\sum\limits_{w\in A_2}\gamma_wt^w, \sum\limits_{y\in
A_3}\nu(y)t^y]=\pt(x)\sum\limits_{y\in A_3}\nu(y)t^y.$$ Then
$$[\sum\limits_{z\in A_1}\lambda_zt^z\pt, \sum\limits_{y\in
A_3}\nu(y)t^y]=\pt(x)\sum\limits_{y\in A_3}\nu(y)t^y.\eqno(5.2)$$

So if $z\ne 0$, then $\ll_z=0$. Moreover, (5.2) becomes
$$[\lambda_0\pt, \sum\limits_{y\in
A_3}\nu(y)t^y]=\pt(x)\sum\limits_{y\in A_3}\nu(y)t^y.\eqno(5.3)$$
So $\ll_0\pt(y)=\pt(x)$ for all $y\in A_3$. Since $\pt$ is
nondegenerate, $x=\ll_0y$.

Hence
$$\pi(\pt)=\ll_0\pt+\sum\limits_{z\in A_2}\gamma_zt^z\eqno(5.4)$$
and $$\pi(t^x)=\nu(\ll_0^{-1}x)t^{\ll_0^{-1}x}.$$

So $\pi(I)\subseteq I'= \span_{\mathbb F}\{t^{\ll_0^{-1} x}, x\in
A\}$. But $I=I'$, we deduce that  $\vep=\ll_0^{-1}\in {\cal E}$.

Therefore
$$\pi(t^x)=\nu(\vep x)t^{\vep x}$$
and
$$\pi(\pt)=\vep^{-1}\pt+\sum\limits_{w\in A_2}\gamma_wt^w$$
where $\vep\in {\cal E}$.

Let $$\eta=\prod_{w\in A_2, w\ne0}\exp(-{\vep\gamma_w \over
\pt(w)}\ad(t^w)),$$ then $\eta^{-1}\pi(\pt)=\vep^{-1}\pt+a$, where
$a=\gamma_0\in\mathbb F,\, \vep\in{\cal E}$. \hfill
$\rule[-.23ex]{1.0ex}{2.0ex}$

Let $\chi\in X(A)$, $\vep\in{\cal E}$, $a, b\in{\mathbb F}$,
$c\in{\mathbb F}^*$, then there is a unique linear map
$$\theta=\theta(\chi, \vep, a, b, c): {\cal D}_1\to {\cal D}_1$$
such that
$$\theta(t^x\pt)=\vep^{-1}\chi(x)t^{\vep x}\pt+(b\pt(x)+a)\chi(x)t^{\vep x}, \quad \theta(t^y)=c\chi(y)t^{\vep y}\eqno(5.7)$$ for all $x\in
A$.

It is straightforward to verify that $\theta$ is an automorphism of
Lie algebra ${\cal D}_1$.

\vskip 1mm
\noindent{\bf Lemma 5.3} {\sl
Let $\pi$ be an automorphism of ${\cal D}_1$, then there exist $b\in{\mathbb F}$,
$c\in{\mathbb F}^*$, $\chi\in X(A)$, $\vep\in {\cal E}$, such that}
$\pi=\theta(\chi, \vep, a, b, c).$

\noindent{\bf Proof.} By Lemma 5.2, we can suppose that
$\pi(\pt)=\vep^{-1}\pt+a$ and $\pi(t^x)=\nu(x)t^{\vep x}$ for some
$a\in{\mathbb F},\, \vep\in{\cal E}$ and $\nu(x)\in \mathbb F$. By
the fact that $[t^x\pt, t^y]=\pt(y)t^{x+y}$ we can deduce that
$\pi(t^x\pt)\in L_{\vep x}$. So
$$\pi({L}_x)={L}_{\vep x}.$$

For each $x\in A$, we can suppose that
$$\pi(t^x\pt)=\vep^{-1}\chi(x)t^{\vep x}\pt+\tau(x)t^{\vep x},\eqno(5.8)$$
$$\pi(t^x)=\ll(x)t^{\vep x},\eqno(5.9)$$
where $\chi(x),\,\ll(x)\in{\mathbb F}^*$, $\tau(x)\in{\mathbb F}$.
Moreover, $\chi(0)=1$ and $\tau(0)=a$ since
$\pi(\pt)=\vep^{-1}\pt+a$. We claim that
$\chi(x+y)=\chi(x)\chi(y)$ for all $x, y\in A$. It suffices to
prove this for $x\ne y$. Then by applying $\pi$ to
$$[t^x\pt, t^y\pt]=\pt(y-x)t^{x+y}\pt\eqno(5.10),$$ we obtain
$$ [\vep^{-1}\chi(x)t^{\vep x}\pt+\tau(x)t^{\vep
x},\ \vep^{-1}\chi(y)t^{\vep y}\pt+\tau(y)t^{\vep y}] $$
 $$=\pt(y-x)\big(\vep^{-1}\chi(x+y)t^{\vep(x+y)}\pt+\tau(x+y)t^{\vep
(x+y)}\big). \eqno(5.11)$$ Therefore
$$\chi(x+y)=\chi(x)\chi(y)\eqno(5.12)$$ and
$$\pt(y-x)\tau(x+y)=\pt(y)\tau(y)\chi(x)-\pt(x)\tau(x)\chi(y).\eqno(5.13)$$
Suppose that $\tau(x)=f(x)\chi(x)$, then (5.13) becomes
$$\pt(y-x)f(x+y)=\pt(y)f(y)-\pt(x)f(x).\eqno(5.14)$$
Fix $0\ne x_0\in \F$.  By substituting $x_0$ and $2x_0$ for $y$ in
(5.14) respectively, we obtain
$$\pt(x_0-x)f(x+x_0)=\pt(x_0)f(x_0)-\pt (x)f(x)\eqno(5.15)$$
and
$$\pt(2x_0-x)f(x+2x_0)=2\pt(x_0)f(2x_0)-\pt (x)f(x).\eqno(5.16)$$
By substituting $x_0$ for $x$ and $x+x_0$ for $y$ in (5.14), we
obtain
$$\pt(x)f(x+2x_0)=(\pt (x)+\pt(x_0))f(x+x_0)-\pt(x_0)f(x_0).\eqno(5.17)$$
Since $\vp$ is non-degenerate, by (5.15)-(5.17), we deduce that
$$
f(x)=\frac{f(2x_{0})-f(x_{0})}{\pt (x_{0})}\pt
(x)+2f(x_{0})-f(2x_{0}).$$ So we can assume that $f(x)=b\pt(x)+a$
for some  $b\in{\mathbb F}$.

By applying $\pi$ to
$$[t^x\pt, t^y]=\pt(y)t^{x+y},$$
we obtain
$$[\vep^{-1}\chi(x)t^{\vep x}\pt+\tau(x)t^{\vep x},\ \lambda(y)t^{\vep y}]=\pt(y)\lambda(x+y)t^{\vep (x+y)}.$$
So
$$\chi(x)\lambda(y)\pt(y)t^{\vep (x+y)}=\pt(y)\lambda(x+y)t^{\vep (x+y)}.$$
Therefore
$$\lambda(x+y)\pt(y)=\chi(x)\lambda(y)\pt(y).\eqno(5.18)$$ By
setting $x=-y$ in (5.18), we obtain $\lambda(y)=c\chi(y)$ if
$y\ne0$, where $c=\lambda(0)$. Therefore $\lambda(x)=c\chi(x)$ for
all $x\in A$ since $\chi(0)=1$.  \hfill
$\rule[-.23ex]{1.0ex}{2.0ex}$

By Lemma 5.2 and Lemma 5.3, we have the following result.

\vskip 1mm
\noindent{\bf Theorem 5.4} {\sl
Set $$\frak{aut}\;({\cal D}_1)=\{\theta(\chi, \vep, a, b, c)\mid
\chi\in X(A), \vep\in{\cal E}, a, b\in{\mathbb F}, c\in{\mathbb
F}^*\},$$ then $$\aut({\cal D}_1)\cong \frak r\rtimes
\frak{aut}\;({\cal D}_1).$$}

Let $\theta(\chi_i, \vep_i, a_i, b_i, c_i),\, i=1, 2$ be two
automorphisms of ${{\cal D}_1}$ defined as in (5.7), then
$$\theta(\chi_1, \vep_1, a_1, b_1, c_1)\theta(\chi_2, \vep_2, a_2, b_2,
c_2)$$$$ =\theta((\chi_1\circ\vep_2)\chi_2,\, \vep_1\vep_2,
\,\vep_2^{-1}a_1+c_1a_2,\, b_1+c_1b_2, \,c_1c_2).\eqno(5.19)$$

Moreover $$\theta(\chi, \vep, a, b, c)^{-1}
=\theta(\chi^{-1}\circ\vep^{-1}, \vep^{-1}, -\vep ac^{-1},
-bc^{-1}, c^{-1}).\eqno(5.20)$$

The map $X(A)\to \frak{aut}\;({\cal D}_1)$ defined by
$$\chi\mapsto \theta(\chi, 1, 0, 0, 1)$$
is an injective homomorphism. Let $N$ be the image of this
homomorphism, then $N$ is a normal subgroup of $\frak{aut}\;({\cal
D}_1)$ by (5.19) and (5.20).

Set $${\frak a}=\{\theta(1, 1, a, b, 1)\mid a, b\in{\mathbb F}\},\
{\frak c}=\{\theta(1, 1, 0, 0, c)\mid c\in{\mathbb F}^{*}\}.$$
Clearly
$${\frak a}\cong {\mathbb F}\times {\mathbb F}, \quad {\frak c}\cong {\mathbb F}^*,$$
and ${\frak a},\, {\frak c}$ are subgroups of $\frak{aut}\;({\cal
D}_1)$. Moreover  $N{\frak a}{\frak c}$ is a normal group of
$\frak{aut}\;({\cal D}_1)$.

Denote by ${\mathcal A}$ the image of the endomorphism of
$\frak{aut}\;({\cal D}_1)$ defined by $\theta(\chi, \vep, a, b,
c)\mapsto \theta(1, \vep, 0, 0, 1)$. Its kernel is $N{\frak a}{\frak
c}$ and its restriction to ${\mathcal A}$ is the identity map. Hence
$\frak{aut}\;({\cal D}_1)\cong (N{\frak a}{\frak c})\rtimes
{\mathcal A}$ and ${\cal A}\cong {\cal E}$. Therefore we have the
following result.

\vskip1mm
\noindent{\bf Theorem 5.5}
$$\frak{aut}\;({\cal D}_1)\cong X(A){\frak a}{\frak
c}\rtimes{{\cal E}}.$$

\noindent\bf{\footnotesize Acknowledgements}\quad\rm
 {\footnotesize  Project is supported by the NNSF (Grant 10671027,
10701019), the ZJNSF(Grant Y607136, D7080080),
the "Qianjiang Excellence Project" (No. 2007R10031) and the "New Century 151 Talent Project" (2008) of Zhejiang Province .
Authors are grateful to the referee for
correction in some errors and invaluable suggestions.}\\[4mm]

\noindent{\bbb{References}}
\begin{enumerate}

{\footnotesize \item  Arbarello E, DeConcini C, Kac V G, Procesi C.
Moduli spaces of curves and representation theory, Commun
Math Phys, 1988, 117: 1-36

\item Benkart G M, Moody R V. Derivations, central extensions and affine Lie algebras,
Algebras Groups Geom, 1986, 3(4): 456--492

\item  Billig Y. Representations of the twisted
 Heisenberg-Virasoro algebra at level zero, Canad Math Bull, 2003, 46(4): 529--537

\item  Dokovi\'{c} D \v{Z}, Zhao Kaiming.  Generalized Cartan type $W$
Lie algebras in characteristic zero, J Algebra, 1997, 195: 170-210

\item  Dokovi\'{c}, D \v{Z}, Zhao Kaiming.  Derivations, isomorphisms, and second cohomology
of generalized Witt algebras, Trans Amer Math Soc, 1998, 350(2): 643-664

\item Farnsteiner R. Derivations and extensions of finitely generated graded Lie
algebras, J. Algebra, 1988, 118(1): 34-45

\item  Fabbri M A, Moody R V. Irreducible representations
of Virasoro-toroidal Lie algebras,  Comm Math Phys, 1994, 159(1): 1--13

\item  Fabbri M A, Okoh F.   Representations of
Virasoro-Heisenberg algebras and Virasoro-toroidal algebras,
Canad J Math, 1999, 51(3): 523--545

\item  Jiang Qifen, Jiang Cuipo. Representations of the twisted Heisenberg-Virasoro algebra and
the full Toroidal Lie Algebras, Algebra Colloquium, 2007, 14(1): 117--134

\item  Kawamoto N. Generalizations of Witt algebras over a field of characteristic zero,
Hiroshinma Math J,  1986, 16: 417-426

\item Liu Dong, Jiang Cuipo. Harish-Chandra modules over the twisted Heisenberg-Virasoro algebra,
J Math Phys, 2008, 49(1): 1-13

\item  Lv Rencai, Zhao Kaiming.  Classification of irreducible weight modules over the twisted
 Heisenberg-Virasoro algebra, math-ST/0510194

\item  Pianzola A.   Automorphisms of toroidal Lie algebras and their
central quotients,  J. Algebra Appl, 2002, 1(1): 113--121

\item Shen Ran, Jiang Qifen, Su Yucai.  Verma Modules Over the Generalized Heisenberg-Virasoro Algebra, Comm. Alg, 2008, 38(4): 1464-1473.

\item Su Yucai, Zhao Kaiming. Simple algebras of Weyl type, Science in China A, 2001, 44: 419-426.

\item Su Yucai, Zhao Kaiming. Generalized Virasoro and super-Virasoro algebras and modules
 of the intermediate series,  J Alg, 2002, 252: 1--19

\item Su Yucai, Zhao Kaiming. Structure of Lie algebras of Weyl type,  Comm Alg, 2004,  32: 1052--1059

\item Xue Min, Lin Weiqiang, Tan Shaobin. Central extension,
derivations and automorphism group for Lie algebras arising from the
2-dimensional torus, Journal of Lie Theory, 2005, 16:
139-153  }

\end{enumerate}

\end{document}